\documentclass{amsart}

\usepackage[utf8]{inputenc}
\usepackage[T1]{fontenc}
\usepackage{amsmath,amssymb,amsthm}
\usepackage{graphicx}
\usepackage{hyperref}
\usepackage{url}
\usepackage{enumerate}
\usepackage{subcaption}
\usepackage{float}
\usepackage{qrcode}

\hypersetup{
    colorlinks=true,
    linkcolor=magenta,
    filecolor=magenta,      
    urlcolor=magenta,
    citecolor=magenta,
    pdftitle={ReShape: a Collaborative Art Experience},
    pdfauthor={Hugo Parlier and Bruno Teheux},
    pdfsubject={Art, Mathematics, Collaboration, Sonification},
}

\title{ReShape: a Collaborative Art Experience}

\author{Hugo Parlier}
\address{Department of Mathematics, University of Fribourg}
\email{hugo.parlier@unifr.ch}

\author{Bruno Teheux}
\address{Department of Mathematics, University of Luxembourg}
\email{bruno.teheux@uni.lu}

\date{\today}

\begin{document}

\maketitle

\section{Introduction}

The genesis of this project emerged from observing the remarkable balance achieved by abstract artists in their compositions using minimal geometric elements (including Aurélie Némours, Theo van Doesburg, Vilmos Huszár, Vassily Kandinsky, and Piet Mondrian). Their works demonstrate how fundamental graphic components—squares, circles, and points—combined with restricted color palettes can create compelling artistic expressions.

Intrigued by the intrinsic harmony of these works, we designed a crowdsourcing art initiative to explore this phenomenon collaboratively. Over a twelve-month period at various collection points, visitors created geometric compositions through an interactive tablet interface. The creative process was structured by offering participants a curated selection of elements: curves, geometric shapes, and pixel-based designs.

The collected contributions underwent various aggregation processes, forming the foundation for a unique interdisciplinary collaboration. Three composers transformed these visual data sets into distinct musical compositions through various sonification techniques. Working closely with these composers, we developed corresponding visualizations for the final performance at the \href{https://www.fnr.lu/research-with-impact-fnr-highlight/flashback-sound-of-data-where-science-meets-music/}{\textit{Sound of Data Project}}—a collaborative venture between the \href{http://www.uni.lu}{University of Luxembourg}, the \href{http://www.fnr.lu}{Luxembourg National Research Fund} (FNR), the \href{http://www.list.lu}{Luxembourg Institute of Science and Technology} (LIST), and \href{https://rocklab.lu}{Rocklab}, presented as part of \href{https://esch2022.lu/fr/}{Esch2022 European Capital of Culture}.

This report documents the project's conceptual framework, implementation process, collaborative dynamics, and outcomes, while also identifying potential areas for future investigation.

\section{Data Collection Methodology}

\begin{figure}[ht]
\centering
\includegraphics[width=0.5\textwidth]{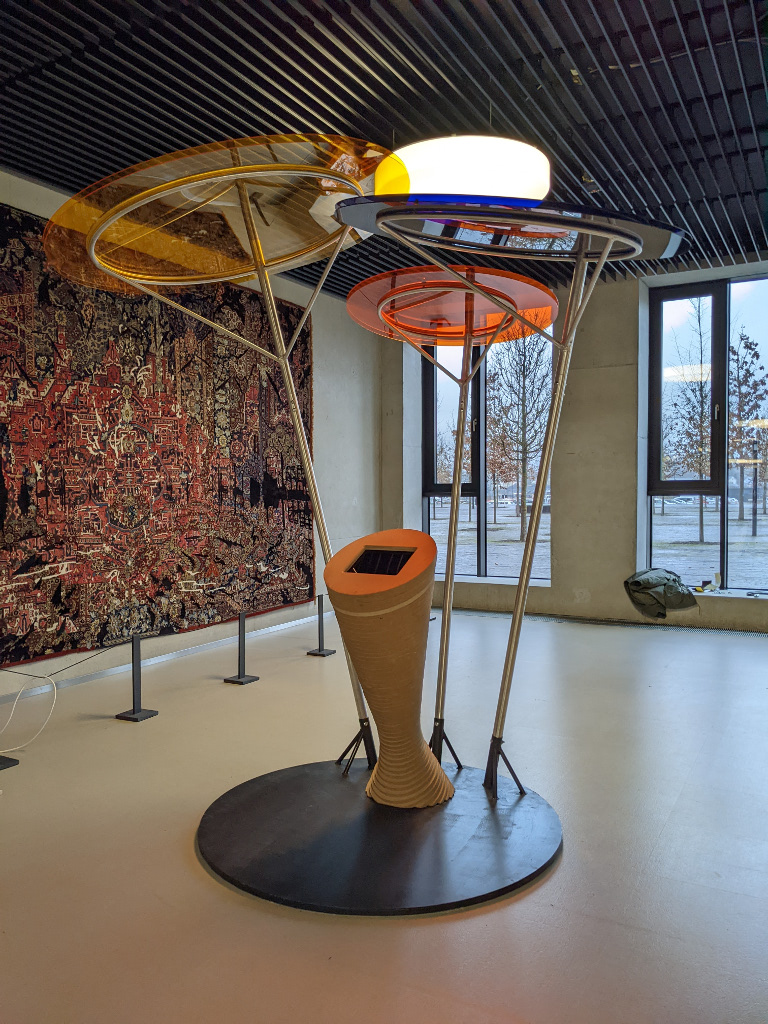}
\caption{Tablet Stand for Data Collection}
\label{fig:tablet_stand}
\end{figure}

The data collection phase was conceived as an outreach activity, designed to engage the public with
mathematical principles through artistic expression. To facilitate this engagement, we collaborated
with a commissioned designer to create an aesthetically appealing tablet stand (Fig.~\ref{fig:tablet_stand}), enabling data collection across multiple locations.

A bespoke application was commissioned, featuring a minimalist interface that encouraged intuitive interaction while maintaining robust backend functionalities for data storage and retrieval. The simplicity of the interface was designed to ensure an engaging user experience.

Data collection sites were strategically selected to reach diverse audiences. These included the \href{http://www.uni.lu}{University of Luxembourg}, \href{https://opderschmelz.lu/}{Centre Culturel Opderschmelz}, \href{https://www.rotondes.lu/}{Rotondes}, \href{https://www.neimenster.lu/}{Neimënster Abbey}, the \href{https://www.science-center.lu/fr}{Luxembourg Science Center}, various local high schools, and the \href{https://www.uni.lu/llc-en/}{Luxembourg Learning Centre}.

The project's reach extended beyond installations in local venues to significant cultural events,
including the Science Meets Music event~\cite{vid:1} and
the Researchers Days 2022~\cite{vid:2}. Most notably, we
conducted a ten-day \href{https://math.uni.lu/outreach/recreate/}{data collection and mediation
session} at the {Luxembourg Pavillon of the World Expo
in Dubai in 2022}, significantly expanding the project's international scope.

\begin{figure}[ht]
\centering
\includegraphics[height=10cm]{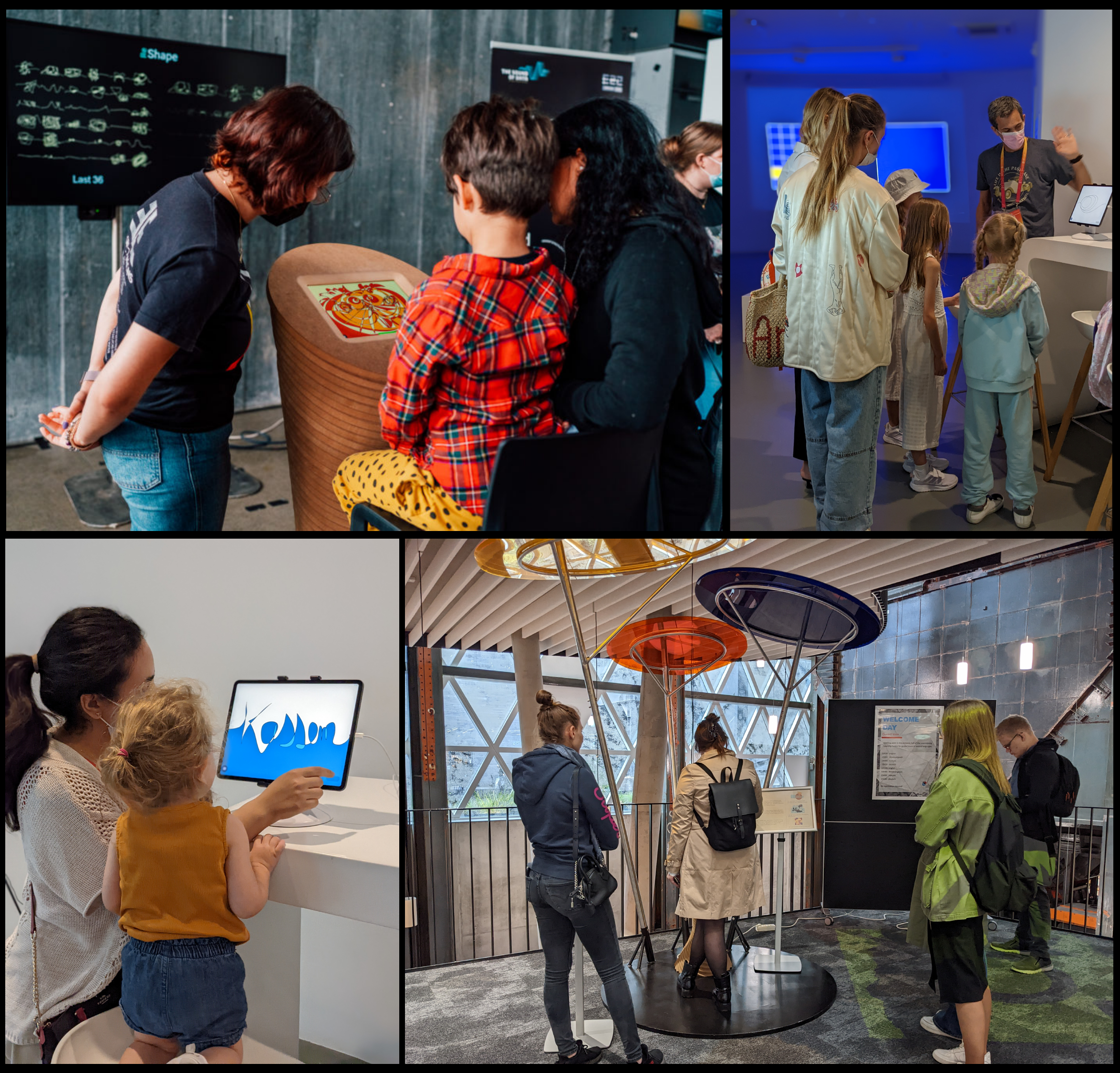}\quad
\caption{Data Collection Events}
\label{fig:data_collection}
\end{figure}

Throughout the data collection process, we prioritized mediated interactions, providing guidance and
explanation to participants (Fig.~\ref{fig:data_collection}; see also \cite{vid:3} for a short illustrative video). This approach created natural opportunities for informal discussions
between researchers and participants about the mathematical principles underlying the artistic
creation process and the various possible aggregation techniques.

By November 2022, the project had amassed a substantial dataset of approximately 20,000 participant-generated drawings, which formed the foundation for the subsequent sonification phase.

\section{Data set}

The data collection interface offered participants three distinct drawing modalities:

\begin{enumerate}
\item \textit{Curve-based compositions}. Participants created drawings by generating curves that traversed the iPad screen from left to right. The regions formed by these curves were automatically colored based on two parameters: the mathematical winding number of the curves and the participant's color selection from a predefined palette. The constraint of left-to-right traversality ensured consistency in the data structure while allowing for significant creative variation.

\item \textit{Shape compositions.} This modality invited participants to construct abstract compositions using combinations of three basic geometric shapes: circles, triangles, and quadrilaterals. Participants could freely draw these shapes, and the regions formed by their intersections were colored according to a mathematical principle. This approach allowed for exploration of compositional balance through fundamental geometric elements.

\item \textit{Grid-based pixel art.} The third modality consisted of fixed-size grid structures where participants created pixel art compositions. This format provided a discrete approach to geometric composition, contrasting with the continuous nature of the other two modalities.
\end{enumerate}

\begin{figure}[ht]
\centering
\includegraphics[height=4cm]{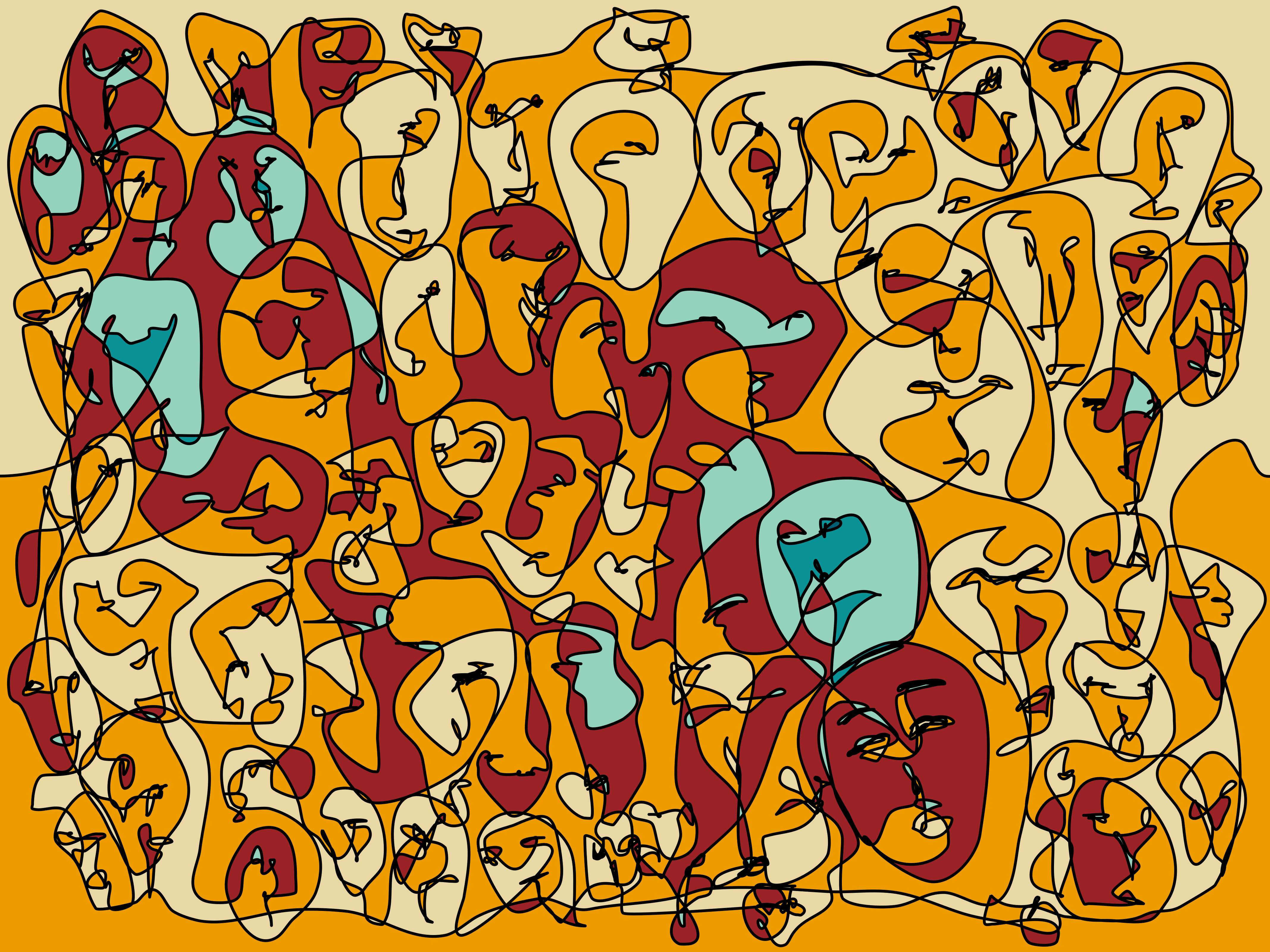}\qquad
\includegraphics[height=4cm]{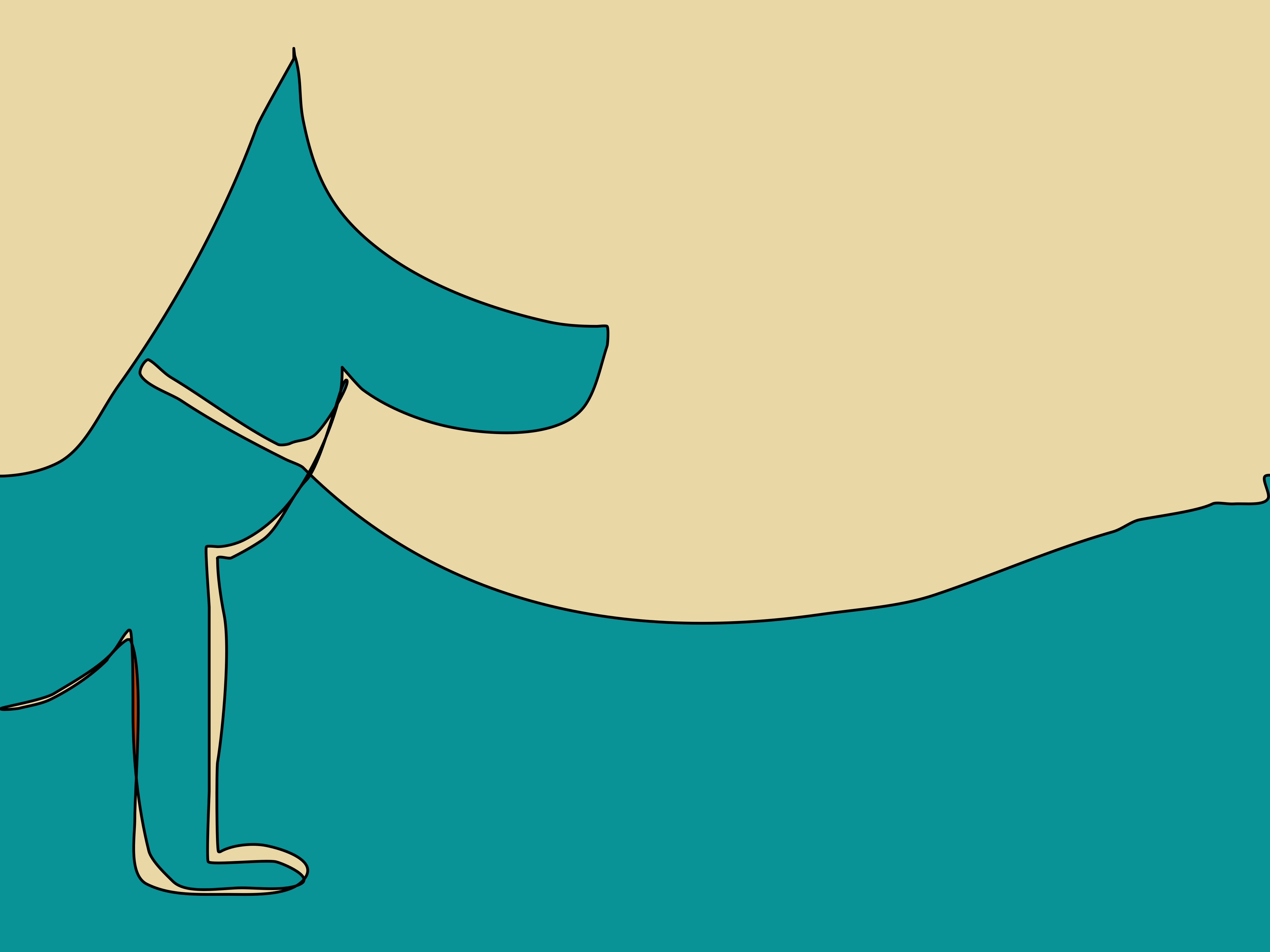}

\includegraphics[height=4cm]{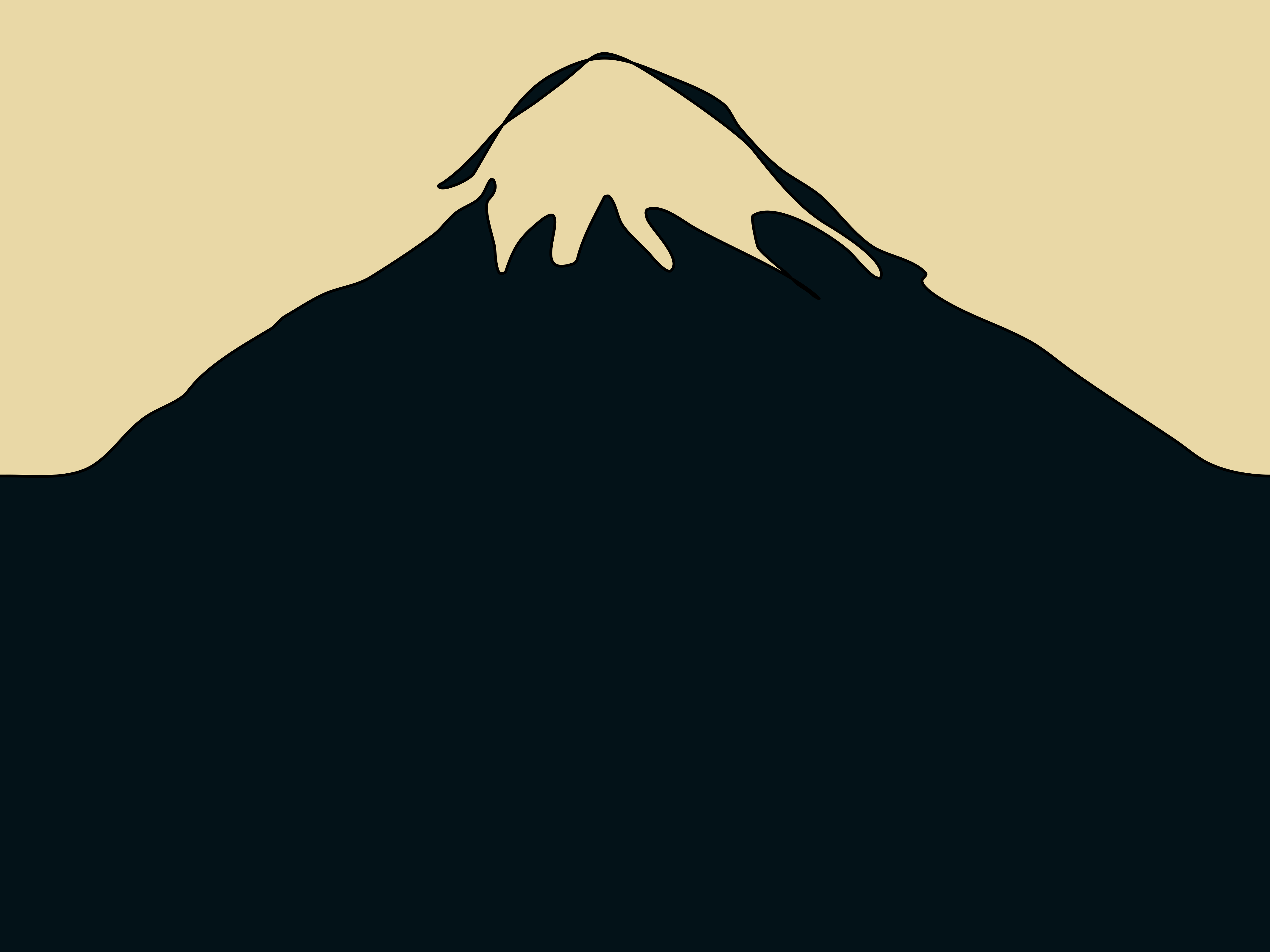}\qquad
\includegraphics[height=4cm]{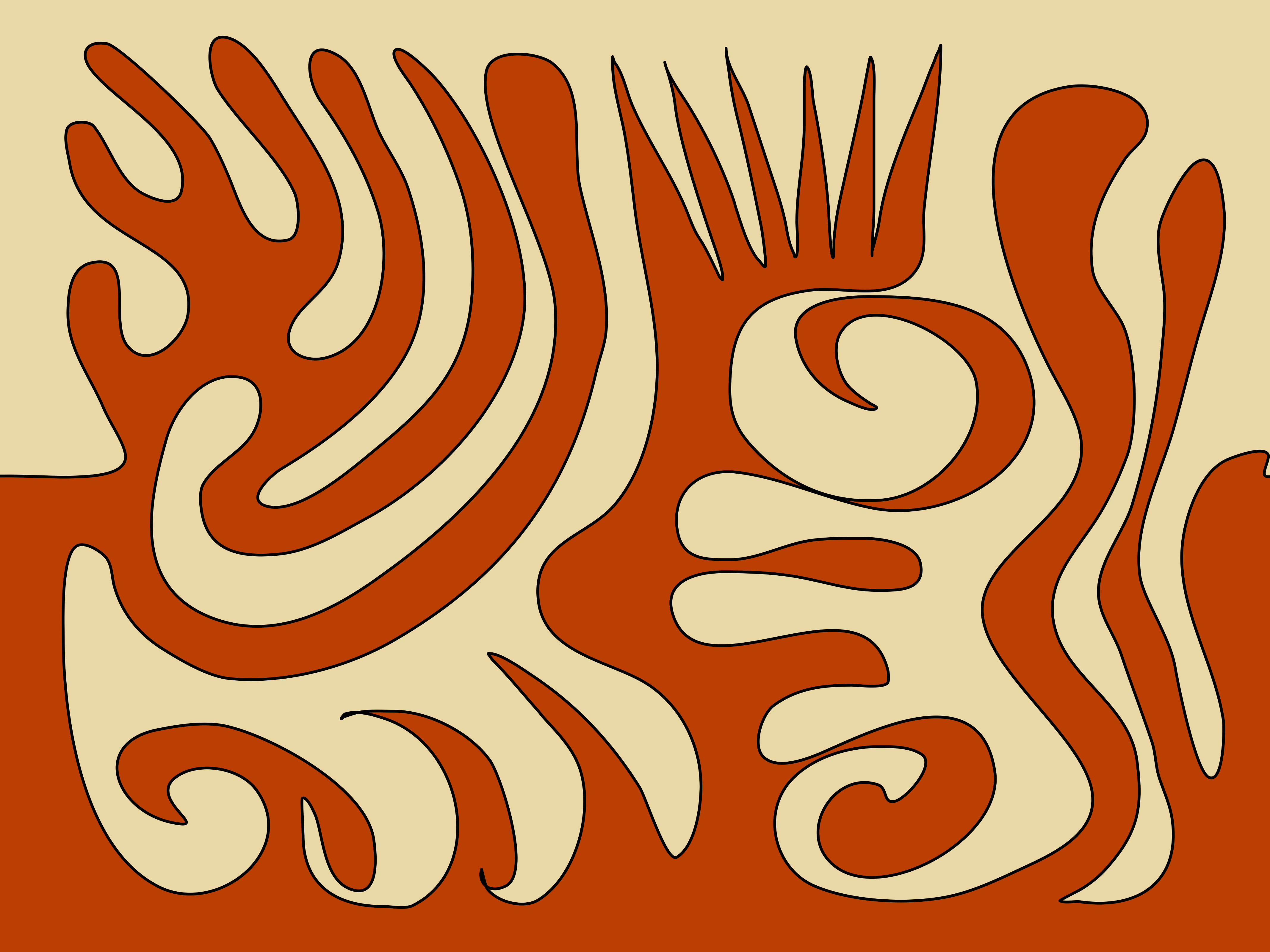}

\includegraphics[height=4cm]{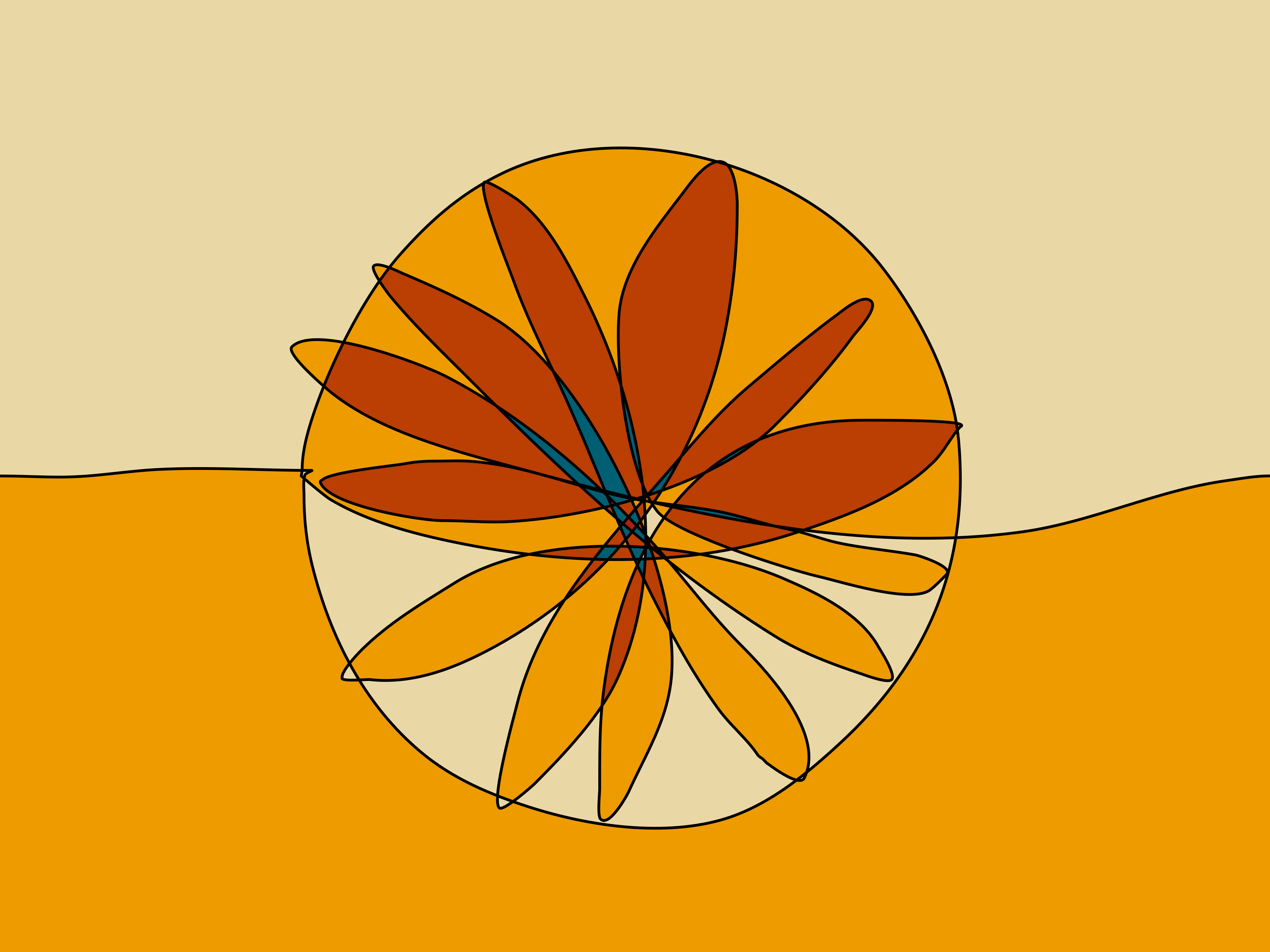}\qquad
\includegraphics[height=4cm]{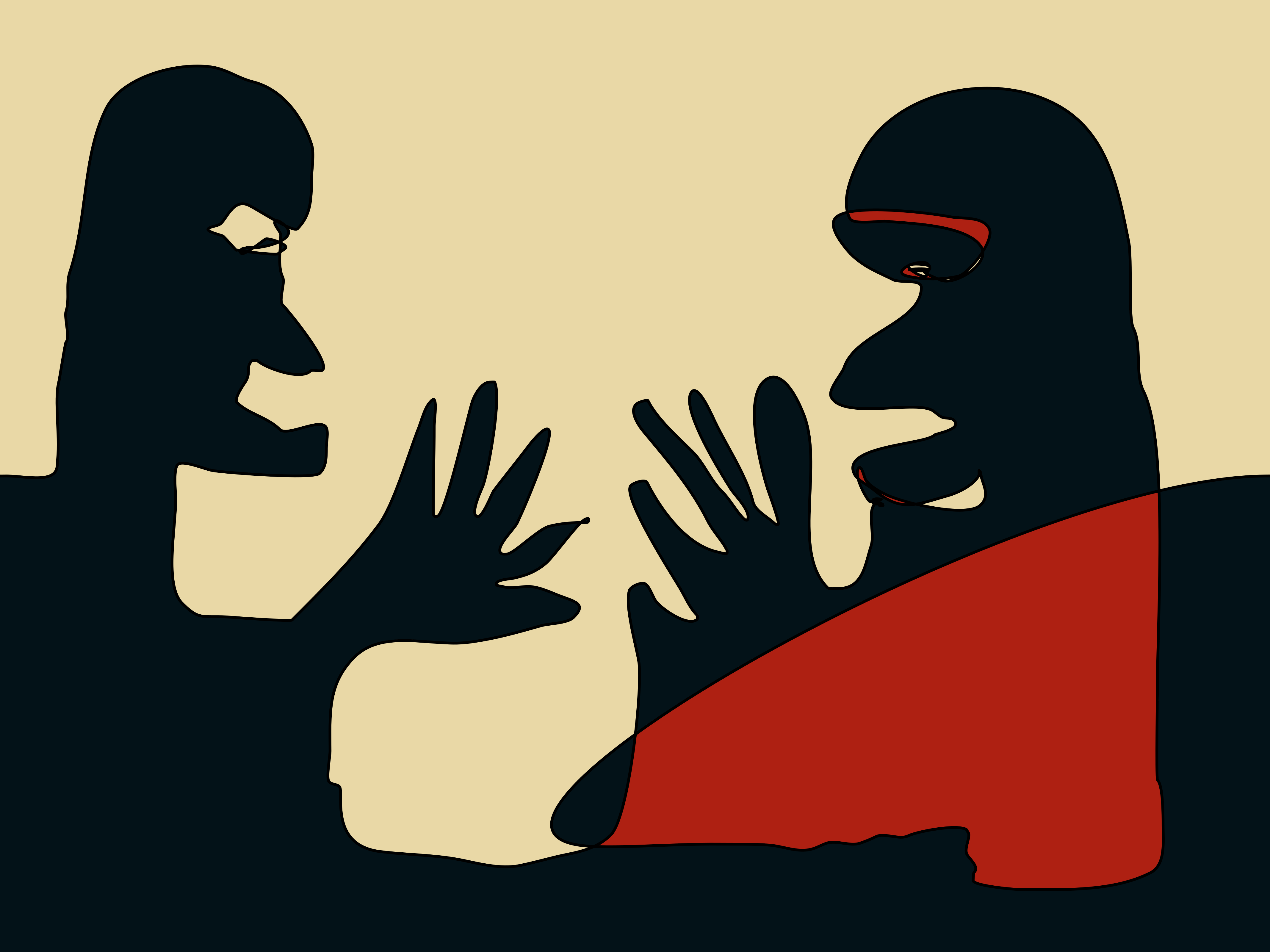}
\caption{Samples of curve-based drawings}
\label{fig:curve_samples}
\end{figure}

\begin{figure}[ht]
\centering
\includegraphics[height=4cm]{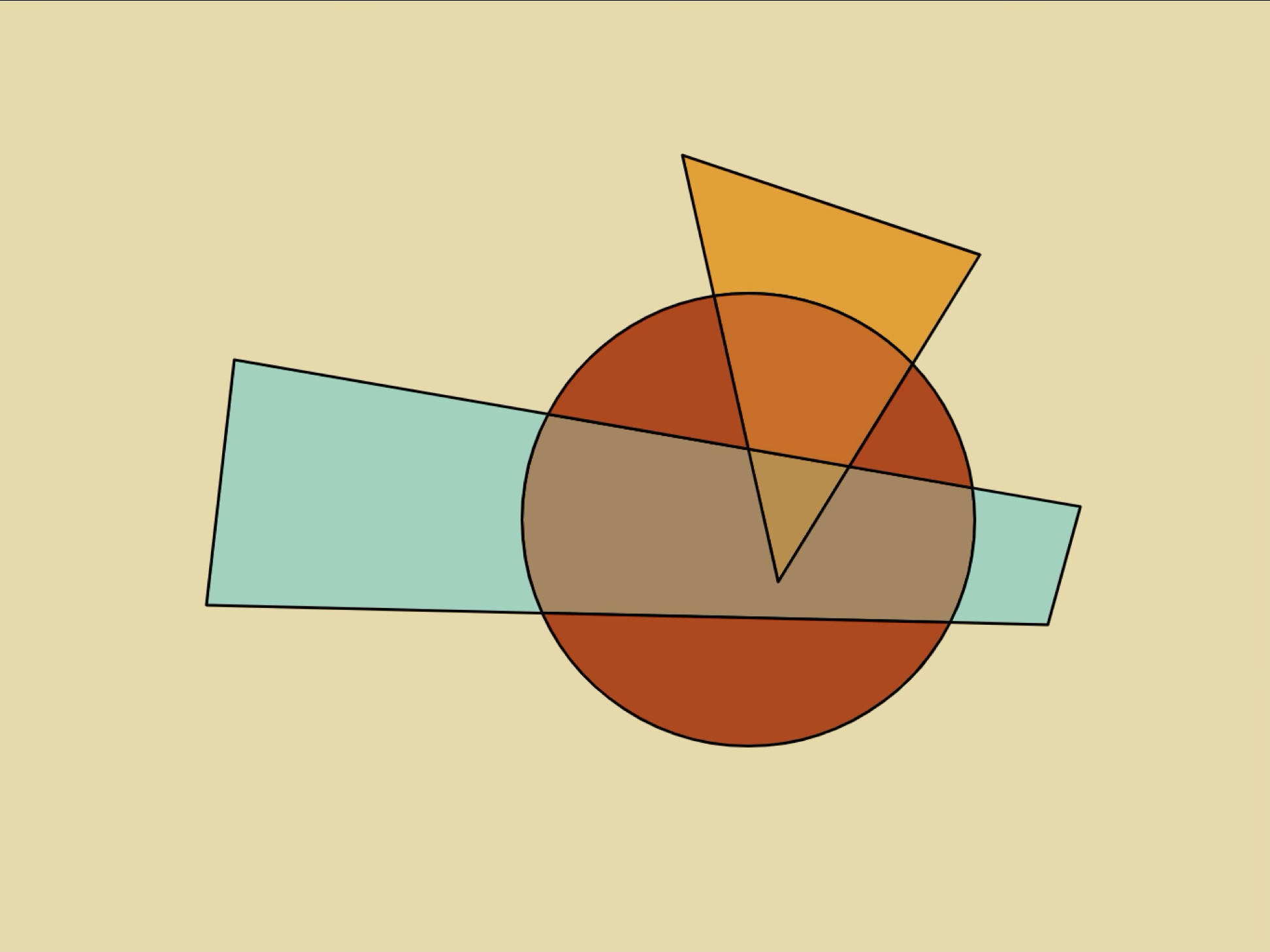}\qquad
\includegraphics[height=4cm]{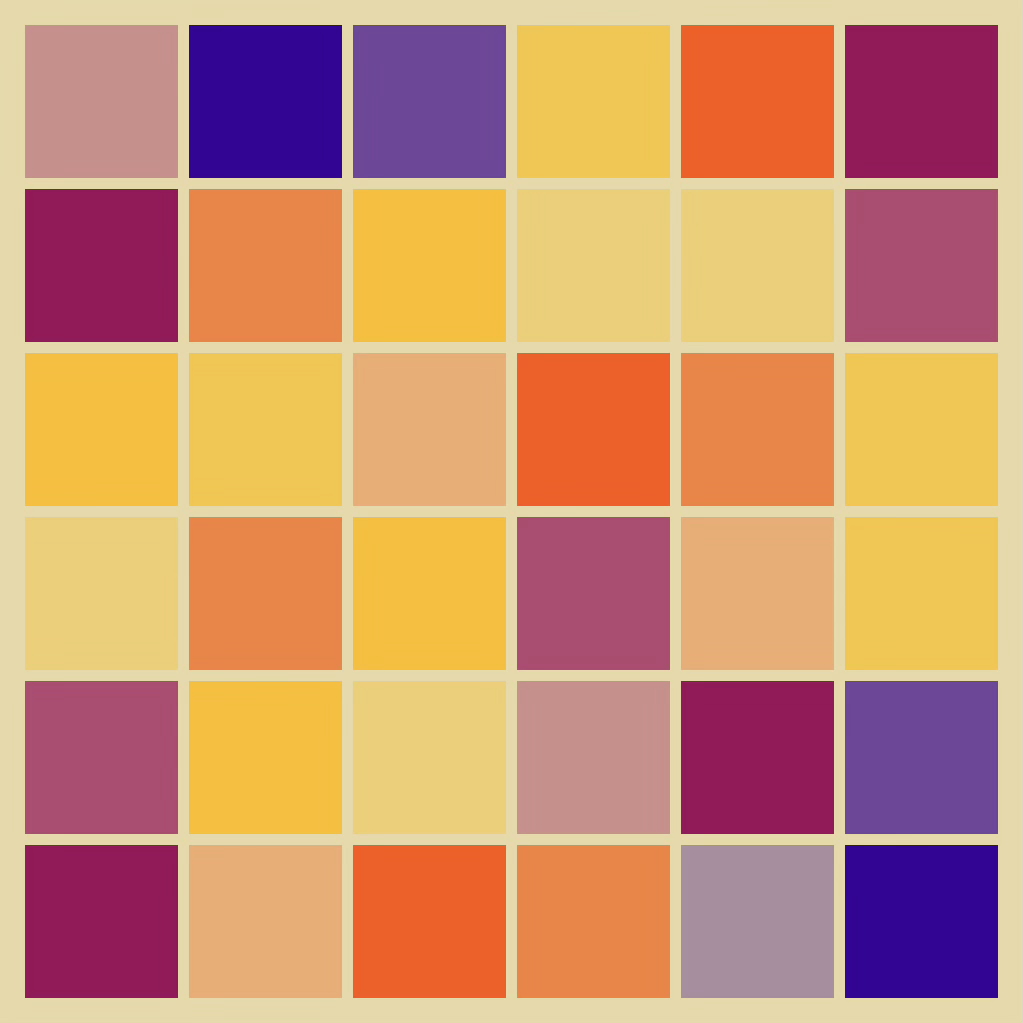}
\caption{Samples of shapes and pixel based drawings}
\label{fig:shape_sample}
\end{figure}

Across all three drawing types, participants worked with a curated and limited color palette. This
constraint not only ensured consistency across the collected datasets but was also inline with the
minimalist principles observed in the geometric abstract art that inspired the project. See
Fig.~\ref{fig:curve_samples} for samples of collected curved-based drawings, and Fig.~\ref{fig:shape_sample} for shapes-based and pixel based drawings.

\section{Aggregation}

The transformation of our collected drawings into collective visual expressions formed a crucial bridge between public participation and the artistic creation of the music by composers. We explored various aggregation methods, not as analytical tools, but as means to create new content that would capture the collective creative expression of all participants.

\subsection{Curve-based drawings.} 
\begin{figure}[ht]
\centering
\includegraphics[width=0.8\textwidth]{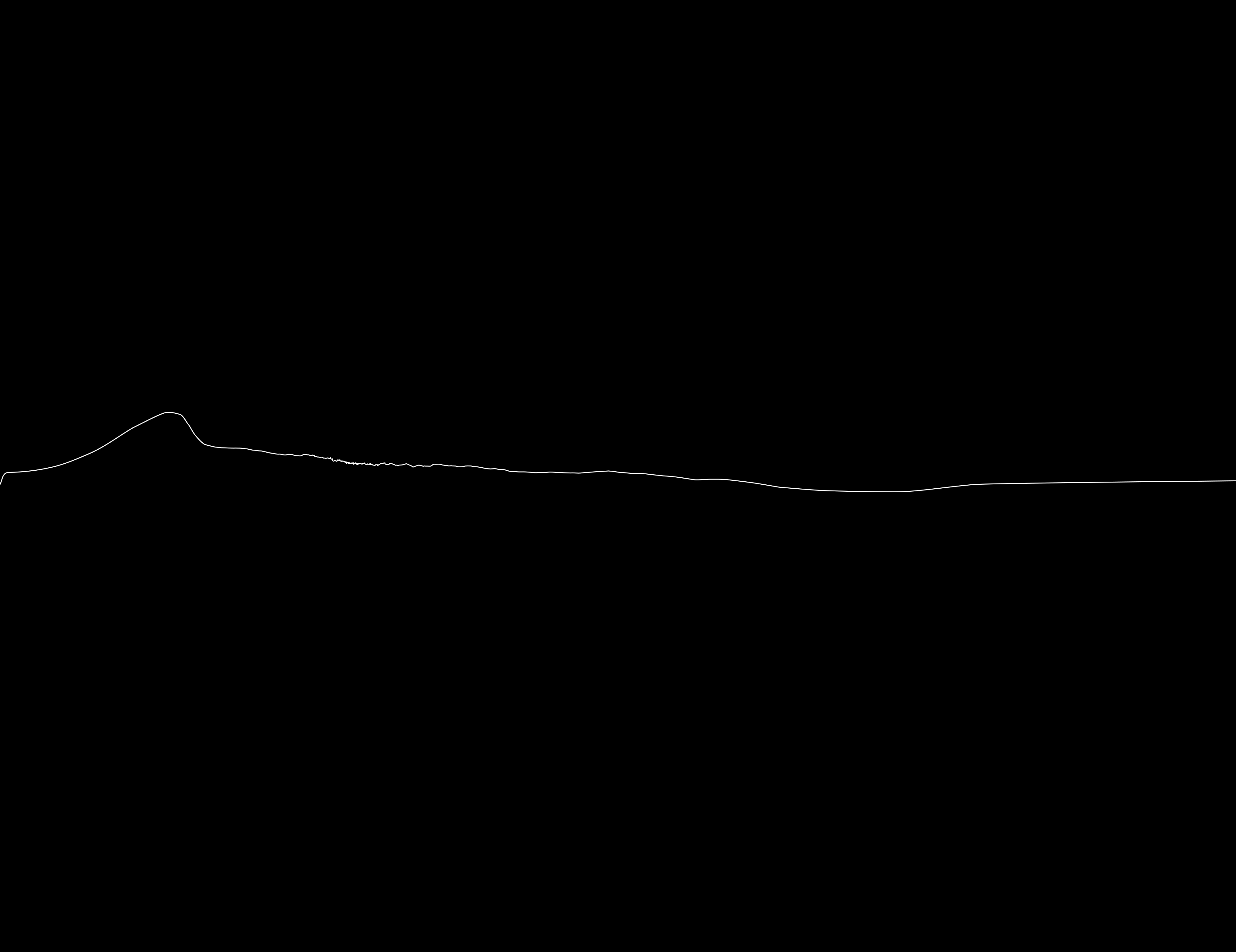}
\caption{Average curve}
\label{fig:average_curve}
\end{figure}

For the curve-based drawings, we developed two complementary approaches. The first method synthesized an 'average curve' by using the mathematical structure of the connected Bézier curves, creating a visual representation that embodied the drawing trajectories of thousands of participants.
Fig.~\ref{fig:average_curve} depicts the computed average, and a video showing a zoom in a small
portion of the curve can be watched at \cite{vid:4}. 

Our second approach explored the space between curves by implementing the so-called Fréchet distance
as a measure, allowing us to visualize shortest paths between different drawings. As a side-note, Fréchet distance is hot research topic where in particular mathematicians are interested in computing it efficiently \cite{buchin2017, cheng2025}. Despite the
computational complexity, we succeeded in approximating this metric space, revealing fascinating
visual journeys through the space of all collected curves. A video showing an animation of a
shortest path in this space of curves can be watched at
\cite{vid:5}.

\subsection{Shape-based drawings.} 
\begin{figure}[H]
\centering
\includegraphics[width=0.8\textwidth]{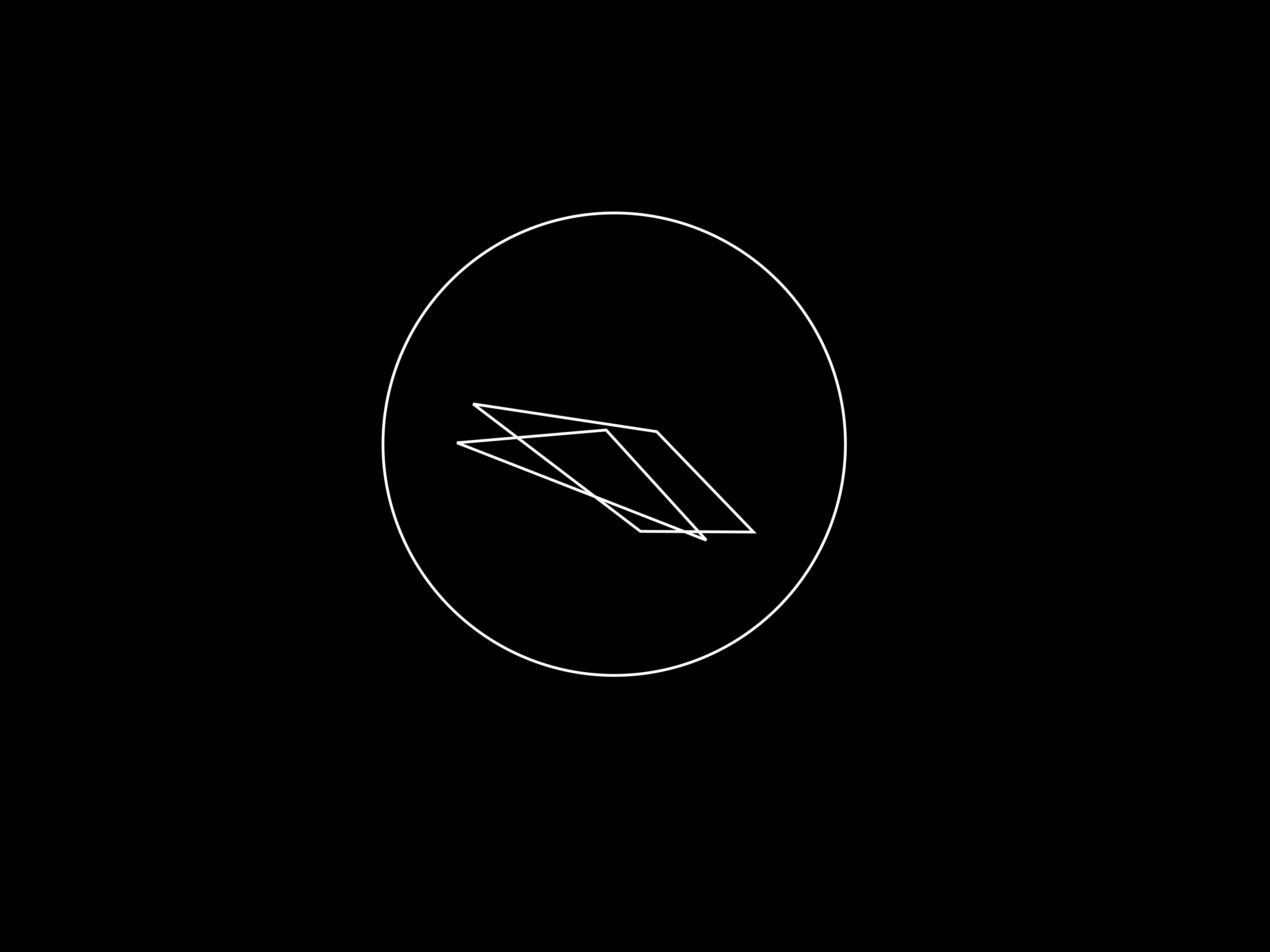}
\caption{Average shapes}
\label{fig:average_shapes}
\end{figure}

We aggregated the shape-based compositions by computing an
average shape for each type of proposed geometry (see Fig.~\ref{fig:average_shapes}). This process transformed the individual circles,
triangles, and quadrilaterals into collective shapes that carried traces of all contributions,
creating new geometric forms at their barycenter.

\textbf{Pixel art.} For the pixel art compositions, we chose to preserve the discrete nature of the
drawing experience through dynamic visualization. By generating a time-lapse of heat maps of the
public action on the pixels, we created a flowing visual narrative that revealed how thousands of
participants collectively activated the pixel grid over time. This approach produced an evolving
artwork that celebrated the cumulative creative process. See
\cite{vid:6} for a video of an animated heat map of
the public action on pixels.

\section{Musical Interpretation, Sonification and Visualization}

The visual data collected through public participation served as the foundation for three distinct musical compositions. In this phase of the project, we adopted a deliberately open approach, entrusting the composers with complete artistic freedom in their interpretation of the geometric datasets. Each composer developed their own unique vision for translating the visual elements into musical compositions.

We worked closely with each composer to identify which aspects of the geometric data resonated most strongly with their artistic vision, and then developed the appropriate data extraction methods to support their creative process.

In parallel to the musical composition process, we developed corresponding visual representations. These visualizations were designed to maintain coherence between the original data, its musical interpretation, and its final presentation, creating a complete audiovisual experience that honored both the mathematical foundations of the project and the composers' artistic interpretations.

\subsection{Mean curve and collective pixel beat}

\href{https://www.hester-1.com/}{Hester-1}'s composition {``Mean curve and collective pixel beat''}
 integrates both geometric and temporal aspects of the collected data. The composition unfolds in two distinct movements, each exploring different dimensions of the collective creative process (see \cite{vid:7}).
\begin{enumerate}
\item \textbf{First Movement: Polyphonic Curves and Human Generated Chaos.} The first movement translates the mathematical properties of the mean curve—its spatial coordinates and tangent data—into a polyphonic arrangement for bowed strings. This foundational layer is enriched by incorporating sonified data from ten randomly selected individual curve drawings, creating a dialogue between the collective average and individual expressions. The accompanying visualization presents an exploration centered on the mean curve phenomenon. Through a carefully crafted zoom into the averaged curve of all participants' drawings, the visualization reveals two concurrent mathematical realities: the emergent collective behavior manifested in the global average pattern, and the intricate local chaos generated by the mathematical averaging process.
\item 
\textbf{Second Movement: Pixel Rhythms and Collective Heat Map.} The second movement transitions to a rhythmic structure derived from the temporal data of pixel drawing interactions. Each screen tap from the original pixel drawing sessions is transformed into a distinct sound, combining computer-generated tones and violin pizzicato to create what the composer terms a ``collective beat.'' This part of the composition is visually represented through an evolving heat map animation that captures the temporal and spatial distribution of participants' interactions, highlighting areas of concentrated activity.
\end{enumerate}

\subsection{Sounding Shadows ReShaped}

\href{http://www.petermichaelvondernahmer.com/}{Mike von der Nahmer's} ``Sounding Shadows ReShaped''
offers a humanistic interpretation of the curve dataset, prioritizing narrative and emotional
resonance over purely algorithmic sonification (see
\cite{vid:8}). The composition unfolds in three distinct movements, each exploring different aspects of human expression through the collected curves.
\begin{enumerate}
\item 
\textbf{First Movement: The Human Heartbeat.} The first movement draws inspiration from curves reminiscent of cardiac rhythms, translating these natural patterns into a musical heartbeat accompanied by morphing curve visualizations. This biological metaphor creates an immediate connection between the mathematical forms and the human experience.
\item
\textbf{Second Movement: The Feedback Loop.} The second movement emerged through an innovative feedback loop: participants created new curves while listening to von der Nahmer's previous compositions, which were then visualized as fading trace animations. These emotionally-informed drawings served as the basis for a musical re-interpretation, effectively giving voice to the participants' visual responses to his original work. This recursive process created a unique dialogue between composer and the participants, and the visualization of their interactions.
\item
\textbf{Third Movement: Mathematical Tribute.} The final movement pays tribute to selected representative drawings from the dataset, accompanied by visualizations that gradually reveal the mathematical complexity of the coloring process, culminating in an animation of a selected shortest path in the space of curves. This movement represents a synthesis of the project's mathematical and artistic elements.
\end{enumerate}

Throughout the piece, von der Nahmer acknowledges the project's international scope and human-centered approach, informally paying tribute to the diverse collection sites including the World Expo in Dubai. The composition reflects both the mathematical structure and the global, collaborative spirit of the project, creating a musical narrative that cross cultural boundaries.

\subsection{Pixelised Harmony}
\href{https://www.timelord.lu/}{Timelord}'s ``Pixelised Harmony'' (see \cite{vid:9}) is structured around a continuous ambient foundation overlaid with a metamorphic main composition that evolves from chaos to order.

The foundational layer consists of a soundscape derived from the pixel art dataset. Using the six-by-six grid drawings as source material, each pixel interaction was assigned a unique sonic representation. These individual sound elements were then layered to create an ambient texture that runs throughout the piece, providing a sonic backdrop for the main compositional elements.

The primary voice of the composition, derived from the shape-based drawings, unfolds in two distinct phases:
\begin{enumerate}
\item 
\textbf{First Movement: Chaos Unveiled.} The first phase sonifies 120 individual shape-based drawings, with each geometric form (circle, triangle, and quadrilateral) assigned a distinct sound signature. The resulting composition manifests as controlled chaos, directly reflecting the diverse and unpredictable nature of individual contributions.
\item 
\textbf{Second Movement: Emergence of Order.} The second phase implements a moving average algorithm to gradually reveal order within the apparent chaos. Beginning with pairs of drawings, the averaging process progressively incorporates more compositions, sonically representing this convergence through a gradual transformation of disparate notes into a coherent central harmony. This musical evolution mirrors the visual averaging process, where individual geometric compositions merge into collective forms, culminating in the general averages of all the collected shapes.
\end{enumerate}

The visual design parallels the sonic architecture through three integrated layers that work in concert to create a mesmerizing experience. Like a living canvas, a dynamic backdrop animates the pixel art drawings that form the foundational soundscape, creating a constantly evolving visual texture that pulses and shifts with the ambient sound base.

In the foreground, during the chaotic first movement, three distinct families of shapes perform their own choreographed sequence. Like the diverse crowd on a dance floor, each geometric family contributes to the collective energy of the resulting ballet.

As the piece transitions into its second phase, these three parallel dances give way to a hypnotic visualization of the averaging algorithm at work. Within each geometric family, the shapes begin to converge, their individual characteristics gradually dissolving to the average ones.

Through this dynamic interpretation, ``Pixelised Harmony'' demonstrates how elementary mathematical concepts—primarily averaging and convergence—can lead to structural harmony from apparent disorder. This fusion of mathematical objects and concepts with contemporary electronic music illustrates how mathematical principles can transcend their abstract nature to create tribal, engaging experiences that resonate with diverse audiences, and pay a primal and visceral tribute to the universality of mathematics.

\section{Final Performance}

The ReShape project's final presentation formed part of the
\href{https://www.fnr.lu/research-with-impact-fnr-highlight/flashback-sound-of-data-where-science-meets-music/}{Sound
of Data Project} showcase on December 2, 2022, an event of \href{https://esch2022.lu/fr/}{Esch2022
European Capital of Culture}. The showcase took place at \href{https://rockhal.lu/}{Rockhal},
Luxembourg's premier concert venue in Esch-sur-Alzette, drawing an audience of over 500 attendees.
The three compositions and their corresponding visualizations were integrated into a program of
\href{https://www.fnr.lu/research-with-impact-fnr-highlight/flashback-sound-of-data-where-science-meets-music/}{eleven
pieces}, all exploring different approaches to data sonification and visualization. The program was
structured in two parts. The first segment presented all eleven Sound of Data pieces, including our
three ReShape works. This concatenation of different approaches to data sonification provided the
audience with a glimpse of how various types of data can be transformed into artistic experiences.
Our project's focus on geometric and mathematical transformations of public contributions offered a
distinct perspective within this diverse showcase. See
\cite{vid:10} for a after movie of the Sound of Data
closing event.

The second half of the evening featured electronic musician and science trained \href{https://maxcooper.net/}{Max Cooper}, whose selection as the headlining artist provided a fitting contextual framework for the entire Sound of Data Project.

\section{Takeaways}

The ReShape project revealed several key insights about public engagement with mathematics through artistic expression. We identified multiple factors contributing to the project's success in engaging diverse audiences.\\

\noindent\textbf{Synergistic integration of art and mathematics.} The project's approach of combining mathematical and artistic creativity proved particularly effective. While traditional outreach activities often struggle to engage participants in either pure mathematical exercises or standalone artistic contributions, our integrated approach successfully bridged this gap. By framing mathematical concepts within a creative context, we observed significantly higher engagement levels than typically seen in single-domain activities.\\

\noindent\textbf{Role of interactive feedback.} The immediate visual feedback provided through the algorithmic coloring process served as a crucial engagement mechanism. Participants were not merely creating static drawings but engaging with a dynamic system that responded to their input. This interactive element encouraged experimentation with the mathematical constraints and motivated participants to create multiple compositions. The process generated natural curiosity about the underlying mechanisms while creating an engaging feedback loop between artistic intent and mathematical transformation.\\

\noindent \textbf{Human-centered mediation.} The mediated nature of the data collection process proved essential to the project's success. Rather than relying on automated or self-guided interactions, the presence of facilitators maintained focus on the profound human component of mathematical activity. As for mathematics or art, human interactions were an essential ingredient of the drawings collection process.\\

\noindent \textbf{Reframing mathematical engagement.} The project demonstrated an alternative approach to mathematical outreach. Instead of presenting mathematics primarily as a tool or ``the language of science''—requiring explanation and instruction—ReShape positioned mathematics as a creative medium. This shift from explaining mathematics to experiencing it through creative expression resulted in more natural and engaging interactions with mathematical concepts, and foster an update in the audience's perspective on mathematics.\\

\noindent \textbf{Universal accessibility.} The project's success in engaging diverse audiences across generations and cultures highlighted the universal appeal of creative geometry based expression. The combination of basic geometric elements, intuitive drawing interfaces, immediate visual feedback, and a collaborative framework created an accessible entry point for mathematical engagement regardless of participants' background or prior mathematical experience.\\

\noindent \textbf{Creative parallels.} The project's collaborative nature effectively highlighted natural parallels between artistic and mathematical creation. The interdependence of artistic works within broader cultural contexts mirrored the interconnected nature of mathematical theorems. The creative process in both domains emerged as inherently collaborative and cumulative, while the need for creativity in art naturally extended to understanding mathematics as a creative process.

\section{Perspectives}

Our experience with the ReShape project has opened several promising avenues for future investigation and development. The collected data and observations from the project implementation suggest multiple directions for both research and outreach activities. In first instance, the relationship between our projects and inspirational work such as Kandinski's "Point and line to plane"  \cite{kandinsky} or more modern approaches to using crowd sourced as an inspiration for design such as in \cite{koyama2018}, remain largely unexplored.

Initial observations indicate potential cultural signatures in the geometric compositions, particularly when comparing drawings collected in Dubai with those from the so-called "greater region" (Luxembourg and neighboring regions). While this pattern is intriguing, it requires rigorous scientific investigation to validate and quantify these apparent cultural influences on geometric expression. 

The interaction between music and geometric creation presents another fascinating research direction. During sessions where participants created drawings while listening to Mike von der Nahmer's compositions, we observed apparent correlations between specific musical pieces and drawing characteristics. A systematic study of these relationships could reveal interesting connections between auditory stimuli and geometric expression, potentially opening new pathways for understanding cross-modal creative processes.

The project's framing as an artistic collaborative endeavor, rather than a mathematical exercise, appeared to significantly impact public engagement. This observation suggests the need for a longitudinal study examining how such framing affects both immediate participation and long-term attitudes toward mathematics. Understanding these dynamics could inform future approaches to mathematical outreach and education.

Our experience also suggests a promising pedagogical sequence: introducing mathematical concepts through creative engagement before formal explanation. This approach, where meaning precedes technical exposition, shows potential for deeper engagement with mathematical tools when they are subsequently presented in dedicated workshops. This sequencing warrants further investigation as a model for mathematical outreach.

The project has generated a substantial dataset of coherent, human-produced geometric drawings that remains largely unexplored. This collection presents numerous opportunities for future research, including investigations into classification methods, complexity measures, and proximity metrics. Such analyses could provide insights into patterns of human geometric thinking and creativity.

Looking forward, we plan to refine the drawing activity by further simplifying its implementation. This refinement aims to distill the experience to its essential creative and mathematical components, potentially making it more accessible while maintaining its effectiveness as a bridge between artistic and mathematical thinking.

These perspectives suggest that the ReShape project, beyond its immediate success as an outreach initiative, has laid groundwork for multiple strands of future research and development in the intersection of mathematics, art, and public engagement.

\begingroup
\def\refname{Links}

\endgroup

\end{document}